\documentclass[letterpaper,11pt]{article}
\usepackage{color}

\usepackage[small]{titlesec}
\usepackage{theorem,amsmath,amssymb,amscd}
\usepackage[all]{xy}
\usepackage[hypertex]{hyperref}

\theoremstyle{change}
\allowdisplaybreaks
\newcommand{\A}{{\mathbb A}}
\newcommand{\Q}{{\mathbb Q}}

\newcommand{\R}{{\mathbb R}}
\newcommand{\C}{{\mathbb C}}

\newcommand{\p}{\mathfrak p}
\newcommand{\OF}{{\mathfrak o}}
\newcommand{\GL}{{\rm GL}}

\newcommand{\SL}{{\rm SL}}
\newcommand{\SO}{{\rm SO}}

\newcommand{\GSp}{{\rm GSp}}

\newcommand{\SSp}{{\rm Sp}}

\newcommand{\qed}{\hspace*{\fill}\rule{1ex}{1ex}}
\newcommand{\forget}[1]{}
\newcommand{\nl}{

\vspace{2ex}}

\def\qdots{\mathinner{\mkern1mu\raise0pt\vbox{\kern7pt\hbox{.}}\mkern2mu
\raise3.4pt\hbox{.}\mkern2mu\raise7pt\hbox{.}\mkern1mu}}

\newtheorem{lemma}{Lemma.}[section]
\newtheorem{theorem}[lemma]{Theorem.}
\newtheorem{corollary}[lemma]{Corollary.}
\newtheorem{proposition}[lemma]{Proposition.}

\begin{document}
\thispagestyle{empty}

\begin{center}
 {\bf\Large Irreducibility Criteria for Local and Global Representations}

 \vspace{2ex}
 Hiro-aki Narita, Ameya Pitale, Ralf Schmidt

 \vspace{3ex}
 \begin{minipage}{80ex}
  \small{\sc Abstract.}
  It is proved that certain types of modular cusp forms generate irreducible automorphic representation of the underlying algebraic group. Analogous archimedean and non-archimedean local statements are also given.
 \end{minipage}

\end{center}


\vspace{4ex}
{\large\bf Introduction}
\nl
One of the motivations for this note was to show that full level cuspidal Siegel eigenforms generate irreducible, automorphic representations of the adelic symplectic similitude group. Such a result is well known for the case of classical elliptic modular forms.
Given an elliptic cusp form $f$ (of some weight and some level), an adelic function $\Phi_f$ can be constructed (see \cite{Ge}, \S5), which is a cuspidal automorphic form on the adelic group $\GL(2,\A)$; here, $\A$ denotes the ring of adeles of $\Q$. Let $V_f$ be the space of automorphic forms generated by all right translates of $\Phi_f$. In this classical situation it turns out that $V_f$ is irreducible precisely when $f$ is an eigenform for the Hecke operators $T_p$ for almost all primes $p$. The proof uses the strong multiplicity one property for cuspidal automorphic representations of $\GL(2)$.
\nl
For most other types of modular forms, strong multiplicity one, or even weak multiplicity one, is not available. The goal of this note is to show that, under certain circumstances, the automorphic representation $V_f$ is still irreducible, even if multiplicity one is not known. Loosely speaking, this is the case whenever $f$ is a holomorphic type of modular form, and is an eigenfunction for \emph{all} Hecke operators. See Corollary \ref{secondglobaltheorem} for a precise statement. We stress that this result does not prove multiplicity one for full level automorphic forms; for example, under our current state of knowledge, it is still conceivable that two holomorphic Siegel cusp forms of degree $n>1$ have the same weight and the same Hecke eigenvalues for all primes $p$, yet are linearly independent.
\nl
While our results are applicable mainly in a reductive setting, we keep the definitions general enough to include certain non-reductive situations, like Jacobi forms. For completeness we also include analogous local archimedean and non-archimedean irreducibility criteria. All of this is well-known to experts, but for lack of a good reference we found it useful to collect the relevant results in this note.

\section{Local non-archimedean theory}
Let $G$ be a group of td-type, as in \cite{Ca}. We fix a left-invariant Haar measure on $G$. The Hecke algebra $\mathcal{H}(G)$ consists of the locally constant, compactly supported functions $f:\:G\rightarrow\C$, with product given by
$$
 (f_1*f_2)(g)=\int\limits_Gf_1(h)f_2(h^{-1}g)\,dh.
$$
Given an open-compact subgroup $K$ of $G$, the symbol $\mathcal{H}(G,K)$ denotes the subalgebra of $\mathcal{H}(G)$ consisting of left and right $K$-invariant functions. If $A$ is a subset of $G$, let $I_A$ be the characteristic function of $A$. Then ${\rm vol}(K)^{-1}I_K$ is an identity element of $\mathcal{H}(G,K)$. As usual, a representation $\pi$ of $G$ on a complex vector space $V$ is called \emph{smooth} if every $v\in V$ is stabilized by some open-compact $K$. In this case $\mathcal{H}(G)$ acts on $V$ via
$$
 \pi(f)v=\int\limits_Gf(g)\pi(g)v\,dg.
$$
The algebra $\mathcal{H}(G,K)$ acts on the space $V^K$ of $K$-invariant vectors. We say that $v\in V^K$ is an eigenvector for $\mathcal{H}(G,K)$ if for all $f\in\mathcal{H}(G,K)$ there exists a scalar $\lambda(f)$ such that $\pi(f)v=\lambda(f)v$.

\begin{lemma}\label{nonarchheckelemma}
 Let $(\pi,V)$ be a smooth representation of $G$. Let $K$ be an open-compact subgroup of $G$. Then, for $v\in V$ and $g\in G$,
 $$
  \int\limits_K\int\limits_K\pi(k_1gk_2)v\,dk_2\,dk_1=
   \frac{{\rm vol}(K)^2}{{\rm vol}(KgK)}\,\pi(I_{KgK})v.
 $$
\end{lemma}
{\bf Proof:} Let $KgK=\bigsqcup g_iK$, a finite disjoint union. Then
\begin{align*}
 \pi(I_{KgK})v&=\sum_i\int\limits_{g_iK}\pi(h)v\,dh
  =\sum_i\int\limits_{K}\pi(g_ik_2)v\,dk_2.
\end{align*}
Applying $\pi(k_1)$ to both sides and integrating over $K$, we obtain
\begin{align*}
 {\rm vol}(K)\pi(I_{KgK})v&=\sum_i\int\limits_K\int\limits_{K}\pi(k_1g_ik_2)v\,dk_2\,dk_1
 =\sum_i\int\limits_K\int\limits_{K}\pi(k_1gk_2)v\,dk_2\,dk_1.
\end{align*}
The assertion follows.\qed

\begin{proposition}\label{nonarchprop}
 Let $(\pi,V)$ be a smooth representation of $G$ such that every $G$-invariant subspace has a $G$-invariant complement\footnote{For example, this is satisfied if $\pi$ is unitarizable.}. Let $K$ be an open-compact subgroup of $G$, and let $v_0\in V^K$ be a non-zero vector with the following properties.
 \begin{enumerate}
  \item The vectors $\pi(g)v_0$, where $g$ runs through $G$, span $V$.
  \item $v_0$ is an eigenvector for $\mathcal{H}(G,K)$.
 \end{enumerate}
 Then the representation $\pi$ is irreducible.
\end{proposition}
{\bf Proof:} We first show that $V^K$ is spanned by $v_0$. Let $v\in V^K$. By hypothesis i), we can write
$$
 v=\sum_{i=1}^nc_i\,\pi(g_i)v_0\qquad\text{for some }c_i\in\C,\;g_i\in G.
$$
Since $v_0$ is $K$-invariant,
$$
 v=\frac1{{\rm vol}(K)}\sum_{i=1}^nc_i\int\limits_K\pi(g_ik_2)v_0\,dk_2.
$$
Applying $\pi(k_1)$ and integrating over $K$, we obtain
$$
 v=\frac1{{\rm vol}(K)^2}\sum_{i=1}^nc_i\int\limits_K\int\limits_K\pi(k_1g_ik_2)v_0\,dk_2\,dk_1.
$$
By Lemma \ref{nonarchheckelemma} and hypothesis ii), the right hand side is a multiple of $v_0$. This proves our claim that $V^K$ is spanned by $v_0$.

\vspace{3ex}
Now let $W$ be a $G$-invariant subspace of $V$, and let $W'$ be a $G$-invariant complement. Write
$$
 v_0=w+w',\qquad w\in W,\;w'\in W'.
$$
Since $v_0$ is $K$-invariant, the same must be true for $w$ and $w'$. By what we already proved, $w$ and $w'$ are both multiples of $v_0$. It follows from the first hypothesis that $V$ is entirely contained in $W$ or in $W'$. This proves that $V$ is irreducible.
\qed

\nl
{\bf Remark:} Let $F$ be a $p$-adic field with ring of integers $\OF$. Let $G=\GL(2,F)$ and $K=\GL(2,\OF)$. Let $V=|\,|^{1/2}\times|\,|^{-1/2}$ be the parabolically induced representation of $G$ which has the Steinberg representation as its unique subrepresentation and the trivial representation as its unique quotient. Then $\dim V^K=1$, and the translates of a non-zero $v_0\in V^K$ span all of $V$, yet $V$ is reducible. This example shows that the hypothesis about invariant complements in Proposition \ref{nonarchprop} is not indispensible.
\section{Local archimedean theory}
Let $G$ be a real Lie group and $K$ a subgroup. We assume that $K$ is of the form center times compact subgroup\footnote{The reason for not assuming that $K$ itself is compact will become apparent once we consider the Jacobi forms example below.}. An equivalence class of irreducible (finite-dimensional) representations of $K$ will be called a $K$-type. With $\mathfrak{g}$ denoting the Lie algebra of $G$, we have the following simple irreducibility criterion.
\begin{lemma}\label{firstarchprop}
 Let $(\pi,V)$ be a $(\mathfrak{g},K)$-module such that every invariant subspace has an invariant complement (for example, this is satisfied if $\pi$ is unitarizable). Let $\tau$ be a $K$-type with the following properties.
 \begin{enumerate}
  \item $V$ is spanned by the vectors $\pi(X)v$, where $X$ runs through $\mathcal{U}(\mathfrak{g}_\C)$ and $v$ runs through the $\tau$-isotypical component of $V$.
  \item $V$ contains the $K$-type $\tau$ exactly once.
 \end{enumerate}
 Then $\pi$ is irreducible.
\end{lemma}
{\bf Proof:} Let $W$ be an invariant subspace of $V$. Let $W'$ be an invariant complement. By ii), exactly one of $W$ or $W'$ contains the $K$-type $\tau$. Then i) implies that either $V=W$ or $V=W'$.\qed
\nl
As above let $\mathfrak{g}$ be the Lie algebra of $G$, and let $\mathfrak{k}$ be the Lie algebra of $K$. Let $\mathfrak{g}_\C$ and $\mathfrak{k}_\C$ be their complexifications.
We will assume that $K$ is connected, so that $K$-types are in one-one correspondence with the equivalence classes of irreducible representations of $\mathfrak{k}_\C$.
We will assume that $\mathfrak{g}_\C$ admits a direct sum decomposition
\begin{equation}\label{liealgebradecompositioneq}
 \mathfrak{g}_\C=\mathfrak{p}_\C^++\mathfrak{k}_\C+\mathfrak{p}_\C^-
\end{equation}
with Lie subalgebras $\mathfrak{p}_\C^+$ and $\mathfrak{p}_\C^-$ such that $[\mathfrak{k}_\C,\mathfrak{p}_\C^\pm]\subset\mathfrak{p}_\C^\pm$. Given a $K$-type $\tau$, we extend $\tau$ to a representation of $\mathfrak{k}_\C+\mathfrak{p_\C^-}$ on which $\mathfrak{p}_\C^-$ acts trivially, and define the $\mathfrak{g}_\C$-module
\begin{equation}\label{universallowestweightmoduleeq}
 L_\tau:=\mathcal{U}(\mathfrak{g}_\C)\otimes_{\mathcal{U}(\mathfrak{k}_\C+\mathfrak{p}_\C^-)}\tau.
\end{equation}

\begin{proposition}\label{secondarchprop}
 Let $G$ and $K$ be as above. Let $(\pi,V)$ be a $(\mathfrak{g},K)$-module such that every invariant subspace has an invariant complement. Let $v_0$ be a vector in $V$ with the following properties.
 \begin{enumerate}
  \item $V$ is spanned by the vectors $\pi(X)v_0$, where $X$ runs through $\mathcal{U}(\mathfrak{g}_\C)$.
  \item The span of the vectors $\pi(k)v_0$, where $k$ runs through $K$, is an irreducible representation $\tau$ of $K$.
  \item The $K$-type $\tau$ occurs in $L_\tau$ exactly once.
  \item $\pi(\p_\C^-)v_0=0$.
 \end{enumerate}
 Then $\pi$ is irreducible.
\end{proposition}
{\bf Proof:} Let $V_\tau$ be the span of the vectors $\pi(k)v_0$, where $k$ runs through $K$. By iv), $\p_\C^-$ annihilates all of $V_\tau$. Therefore, there exists a map
\begin{equation}\label{secondarchpropeq1}
 L_\tau\cong\mathcal{U}(\mathfrak{g}_\C)\otimes_{\mathcal{U}(\mathfrak{k}_\C+\mathfrak{p}_\C^-)}V_\tau
 \longrightarrow V
\end{equation}
given by $X\otimes v\mapsto\pi(X)v$. By i), this map is surjective. Hypothesis iii) implies that $V$ contains the $K$-type $\tau$ only once. Hence, we can apply  Lemma \ref{firstarchprop} to see that $\pi$ is irreducible.\qed
\nl
{\bf Remark:} Hypothesis iii) is satisfied if $G$ is linear, simple and connected, of hermitian type, $K$ is maximal compact, $\mathfrak{p}_\C^+$ (resp.\ $\mathfrak{p}_\C^-$) denotes the space spanned by the root vectors for the non-compact positive (resp.\ negative) roots, and $\tau$ is equivalent to the minimal $K$-type in a holomorphic discrete series representation.
\nl
As an application of Proposition \ref{secondarchprop}, we consider holomorphic Jacobi forms on $\mathcal{H}\times\C$, as in \cite{EZ}. Let $G^J=\SL(2,\R)\ltimes H(\R)$ be the real Jacobi group, where $H$ is a three-dimensional Heisenberg group; see \cite{BS}, Sect.\ 1.1, for the precise definition. Let $K$ be the subgroup $\SO(2)\times Z$, where $Z\cong\R$ is the center of the Heisenberg group. Then, by \cite{BS}, Sect.\ 1.4, the complexified Lie algebra $\mathfrak{g}_\C$ of $G^J$ admits a decomposition of the type (\ref{liealgebradecompositioneq}).
\nl
Now let $f$ be a Jacobi form of weight $k$ and index $m$. We associate to $f$ the function $\phi_f:\:G^J\rightarrow\C$ given by $\phi_f(g)=(f\big|_{k,m}g)(i,0)$, where the slash action is the same as in Definition 4.1.1 of \cite{BS}. Then the holomorphy of $f$ is equivalent to $R(\p_\C^-)\phi_f=0$, where $R$ denotes right translation; see Proposition 4.1.2 of \cite{BS}. Moreover, $\phi$ transforms by scalars under right translation by elements of $K$. Condition iii) of Proposition \ref{secondarchprop} is also satisfied, as is apparent from the description of the lowest weight representation $\pi_{k,m}^+$ of $G^J$ in \cite{BS}, Sect.\ 3.1. Therefore, we obtain the following corollary (which could be generalized to higher degree Jacobi forms).
\begin{corollary}\label{secondarchpropcor}
 Let $f$ be a holomorphic Jacobi form on $\mathcal{H}\times\C$ of weight $k$ and index $m$. Then, under right translation, the associated function $\phi_f:\:G^J\rightarrow\C$ generates an irreducible $(\mathfrak{g},K)$-module.
\end{corollary}
\nl
We also seek an irreducibility criterion for non-connected Lie groups; this is important for certain global applications which we will give below. Hence assume that $G$ is a Lie group with connected component $G_0$ such that $G_0$ has finite index in $G$. Let $K$ be a compact subgroup of $G$ whose connected component $K_0$ has finite index in $K$. If $(\pi,V)$ is a $\mathfrak{g}_\C$-module, and $\sigma\in K$ is an element, let $(\pi^\sigma,V)$ be the $\mathfrak{g}_\C$-module with underlying space $V$ and action given by
\begin{equation}\label{twistedgkmodule}
 \pi^\sigma(X)v=\pi({\rm Ad}(\sigma)X)v,\qquad X\in\mathfrak{g}_\C,\:v\in V.
\end{equation}
Again we assume that $\mathfrak{g}_\C$ admits a decomposition as in (\ref{liealgebradecompositioneq}).

\begin{proposition}\label{thirdarchprop}
 Let $G$ and $K$ be as above. Let $(\pi,V)$ be a $(\mathfrak{g},K)$-module such that every invariant subspace has an invariant complement. Let $v_0$ be a vector in $V$ with the following properties.
 \begin{enumerate}
  \item $V$ is spanned by the vectors $\pi(X)\pi(k)v_0$, where $X$ runs through $\mathcal{U}(\mathfrak{g}_\C)$ and $k$ runs through $K$.
  \item The span of the vectors $\pi(k)v_0$, where $k$ runs through $K_0$, is an irreducible representation $\tau$ of $K_0$.
  \item The $K_0$-type $\tau$ occurs in $L_\tau$ (defined as in (\ref{universallowestweightmoduleeq})) exactly once.
  \item For $\sigma\in K$, $\sigma\notin K_0$, the $\mathfrak{g}_\C$-modules $L_\tau$ and $(L_\tau)^\sigma$ have no $K_0$-type in common.
  \item $\pi(\p_\C^-)v_0=0$.
 \end{enumerate}
 Then $\pi$ is irreducible.
\end{proposition}
{\bf Proof:} By i),
$$
 V=\sum_{\sigma\in K/K_0}\pi(\sigma)\pi(\mathcal{U}(\mathfrak{g}_\C))v_0.
$$
As in the proof of Proposition \ref{secondarchprop}, there is a surjective map $L_\tau\rightarrow\pi(\mathcal{U}(\mathfrak{g}_\C))v_0$. Hence, by iii), the $K_0$-type $\tau$ occurs in $\pi(\mathcal{U}(\mathfrak{g}_\C))v_0$ exactly once. Using iv), it follows that the $K_0$-type $\tau$ occurs in all of $V$ exactly once.
\nl
Now let $W\subset V$ be an invariant $(\mathfrak{g},K)$-submodule, and let $W'$ be an invariant complement. Then exactly one of $W$ or $W'$ contains the $K_0$-type $\tau$, and hence the vector $v_0$. It follows from i) that either $V=W$ or $V=W'$.\qed
\nl
In this section we have kept our hypotheses general in order to accomodate examples like that of Jacobi forms. For background on highest weight modules and irreducibility in the semisimple setting, see Sect.\ 7 of \cite{Ya}.
\section{Global theory}
Let $F$ be a number field, and let $G$ be an algebraic group defined over $F$. Let $Z$ be the center of $G$. Let $\A$ be the ring of adeles of $F$, and let $\A_\infty$ and $\A_f$ be the subrings of archimedean and finite adeles, respectively. Then we have the adelized groups $G(\A)$, $G(\A_\infty)$ and $G(\A_f)$. Let $K_\infty$ be a fixed maximal compact subgroup of $G(\A_\infty)$. For the present purposes, we call a function $\Phi:\:G(\A)\rightarrow\C$ an automorphic form on $G(\A)$ if the following conditions are satisfied.
\begin{enumerate}
 \item $\Phi$ is smooth, i.e., $\Phi$ is right invariant under some open-compact subgroup $K_f$ of $G(\A_f)$, and for fixed $h\in G(\A)$ the function $G(\A_\infty)\ni g\mapsto\Phi(gh)$ is $C^\infty$.
 \item $\Phi$ is $K_\infty$-finite, i.e., the space spanned by the functions $G(\A)\ni g\mapsto\Phi(g\kappa)$, where $\kappa$ runs through $K_\infty$, is finite-dimensional.
 \item $\Phi(\gamma g)=\Phi(g)$ for all $g\in G(\A)$ and $\gamma\in G(F)$.
 \item There exists a unitary character $\chi$ of $Z(\A)$ such that $\Phi(zg)=\chi(z)\Phi(g)$ for all $g\in G(\A)$ and $z\in Z(\A)$.
\end{enumerate}
The usual definition of automorphic forms includes conditions like slow growth, but for our purposes it is irrelevant whether we include such additional conditions or not. The important fact is that the space of automorphic forms is preserved by right translation $R$ by elements of $G(\A_f)$ and of $K_\infty$. It is also preserved by the action of $\mathfrak{g}$, the Lie algebra of $G(\A_\infty)$, given by
$$
 (R(X)\Phi)(g)=\frac d{dt}\Big|_0\Phi(g\exp(tX)).
$$
This action extends to an action of $\mathfrak{g}_\C$ and the universal envelopping algebra $\mathcal{U}(\mathfrak{g}_\C)$. We call any subspace of the space of automorphic forms invariant under right translation by $G(\A_f)$, $K_\infty$ and $\mathfrak{g}$, irreducible or not, an \emph{automorphic representation} of $G(\A)$ (even though it might not be a representation of $G(\A)$ at all). If $\Phi$ is an automorphic form, then the automorphic representation generated by $\Phi$ is by definition the linear span of all functions obtained by the right action of $G(\A_f)$, $\mathfrak{g}$ and $K_\infty$ on $\Phi$. The automorphic form $\Phi$ is called square-integrable if
$$
 \int\limits_{G(F)Z(\A)\backslash G(\A)}|\Phi(g)|^2\,dg<\infty.
$$
Here, the integration is with respect to a Haar measure on $G(\A)$. The space of square-integrable automorphic forms with respect to a fixed unitary character of the center $Z(\A)$ is invariant under right translation and carries an obvious inner product.

\begin{theorem}\label{maintheorem}
 For each finite place $v$, let $K_v$ be an open-compact subgroup of $G(F_v)$; we assume that almost all, but not necessarily all, $K_v$ are maximal compact in $G_v$. Let $\Phi$ be a square-integrable automorphic form on $G(\A)$ satisfying the following properties.
 \begin{enumerate}
  \item The $(\mathfrak{g},K_\infty)$-module spanned by the functions $R(X)R(k)\Phi$, where $X$ runs through $\mathfrak{g}_\C$ and $k$ runs through $K_\infty$, is irreducible.
  \item For each finite place $v$, the function $\Phi$ is right invariant under $K_v$, and is an eigenvector under the action of the local Hecke algebra $\mathcal{H}(G(F_v),K_v)$.
 \end{enumerate}
 Then the automorphic representation of $G(\A)$ generated by $\Phi$ is irreducible.
\end{theorem}
{\bf Proof:} Let $V$ be the automorphic representation generated by $\Phi$. Let $F$ be a non-zero function in $V$ which is right invariant under $K_v$ for all finite $v$. We can write
\begin{equation}\label{maintheoremeq1}
 F=\sum_{i=1}^nR(X_i)R(\kappa_i)R(h_i)\Phi,\qquad X_i\in\mathcal{U}(\mathfrak{g}_\C),\;\kappa_i\in K_\infty,\;h_i\in G(\A_f).
\end{equation}
Since $\Phi$ is right invariant under almost all local maximal compact subgroups, we may assume that there exists a finite set $S$ of finite places such that each $h_i$ has trivial components outside of $S$. For $v\in S$, since $F$ and $\Phi$ are both $K_v$-invariant,
$$
 F=\frac1{{\rm vol}(K_v)^2}
 \sum_{i=1}^nR(X_i)R(\kappa_i)\int\limits_{K_v}\int\limits_{K_v}R(k_1h_ik_2)\Phi
 \,dk_1\,dk_2.
$$
By Lemma \ref{nonarchheckelemma}, the double integration coincides, up to a multiple, with the action of the double coset $K_v(h_i)_vK_v$ on $\Phi$. By hypothesis ii), this action reproduces $\Phi$ up to a multiple. We may therefore assume that all of the $h_i$'s in (\ref{maintheoremeq1}) have trivial $v$-component. Continuing through all places $v\in S$, we may assume that the elements $h_i$ are not present at all, so that
\begin{equation}\label{maintheoremeq2}
 F=\sum_{i=1}^nR(X_i)R(\kappa_i)\Phi,\qquad X_i\in\mathcal{U}(\mathfrak{g}_\C),\;\kappa_i\in K_\infty.
\end{equation}
Hence $F$ lies in the $(\mathfrak{g},K_\infty)$-module generated by $\Phi$, which is irreducible by hypothesis i). It follows that
\begin{equation}\label{maintheoremeq3}
 \Phi=\sum_{i=1}^mR(Y_i)R(\sigma_i)F,\qquad Y_i\in\mathcal{U}(\mathfrak{g}_\C),\;\sigma_i\in K_\infty.
\end{equation}
We proved that whenever $F$ is non-zero and right invariant under all $K_v$, then a relation of the form (\ref{maintheoremeq3}) holds.
\nl
Now let $W$ be an invariant subspace of $V$. Since $V$ consists entirely of square-integrable functions, there exists an invariant complement $W'$ to $W$ inside $V$. Write
$$
 \Phi=F+F',\qquad F\in W,\;F'\in W'.
$$
Then, evidently, $F$ and $F'$ are both invariant under $K_v$ for all finite $v$. By what we proved, $\Phi$ can be expressed in the form (\ref{maintheoremeq3}), if $F$ is non-zero. A similar relation holds if $F'$ is non-zero. It follows that $\Phi$ is contained in $W$ or in $W'$. This proves the asserted irreducibility of $V$.\qed
\nl
It is desirable to have a more practical criterion for the irreducibility expressed in i) of Theorem \ref{maintheorem}. This is provided by the next result. We assume that $G(\A_\infty)$ is of the form $Z'G'$, where $Z'$ lies in the center of $G(\A_\infty)$, and where $G'$ is a Lie group whose connected component $G'_0$ is of finite index in $G'$. We assume that $K_\infty=(Z'\cap K_\infty)K'$ with a maximal compact subgroup of $G'$ whose connected component $K'_0$ has finite index in $K'$. We also assume that $\mathfrak{g}'_\C$, the complexification of the Lie algebra $\mathfrak{g}'$ of $G'$, admits a direct sum decomposition as in (\ref{liealgebradecompositioneq}), with $\mathfrak{k}$ being the Lie algebra of $K'$. For a $K'_0$-type $\tau$, we define
\begin{equation}\label{universallowestweightmoduleeq2}
 L_\tau:=\mathcal{U}(\mathfrak{g}'_\C)\otimes_{\mathcal{U}(\mathfrak{k}_\C+\mathfrak{p}_\C^-)}\tau.
\end{equation}
If $(\pi,V)$ is a $\mathfrak{g}'_\C$-module, and $\sigma\in K'$ is an element, we define the twisted $\mathfrak{g}'_\C$-module $(\pi^\sigma,V)$ as in (\ref{twistedgkmodule}).

\begin{corollary}\label{secondglobaltheorem}
 For each finite place $v$, let $K_v$ be an open-compact subgroup of $G(F_v)$; we assume that almost all, but not necessarily all, $K_v$ are maximal compact in $G_v$. Let $\Phi$ be a square-integrable automorphic form on $G(\A)$ satisfying the following properties.
 \begin{enumerate}
  \item The span of the vectors $\pi(k)\Phi$, where $k$ runs through $K'_0$, is an irreducible representation $\tau$ of $K'_0$.
  \item The $K'_0$-type $\tau$ occurs in $L_\tau$ exactly once.
  \item For $\sigma\in K'$, $\sigma\notin K'_0$, the $\mathfrak{g}_\C$-modules $L_\tau$ and $(L_\tau)^\sigma$ have no $K'_0$-type in common.
  \item $R(\p_\C^-)\Phi=0$.
  \item For each finite place $v$, the function $\Phi$ is right invariant under $K_v$, and is an eigenvector under the action of the local Hecke algebra $\mathcal{H}(G(F_v),K_v)$.
 \end{enumerate}
 Then the automorphic representation of $G(\A)$ generated by $\Phi$ is irreducible.
\end{corollary}
{\bf Proof:} By Theorem \ref{maintheorem}, it suffices to show that the $(\mathfrak{g},K)$-module generated by $\Phi$ is irreducible. Since $\Phi$ transforms under central elements by a scalar, this is equivalent to saying that the $(\mathfrak{g}',K')$-module generated by $\Phi$ is irreducible. But this is immediate from Proposition \ref{thirdarchprop}.\qed
\nl
As an application of Corollary \ref{secondglobaltheorem}, consider a Siegel modular cusp form $F$ of degree $n$ with respect to the full modular group. We allow the vector-valued case, i.e., $F$ takes values in an irreducible, holomorphic (polynomial) representation $\tau$ of $\GL(n,\C)$. As in \cite{AS}, we associate with $F$ an adelic function $\Phi:\:G(\A)\rightarrow\C$, where $G=\GSp(2n)$ and $\A$ is the ring of adeles of $\Q$. Note that $\Phi$ is complex-valued, even if $\tau$ is more than one-dimensional; in this case the function $\Phi$ is not canonical, but the space of right translates of $\Phi$ by elements of $K'_0\cong U(n)$ is, and is isomorphic to $\tau$.

\begin{corollary}\label{secondglobaltheoremcor1}
 Let $F$ be a vector-valued holomorphic Siegel modular cusp form of degree $n$ and full level. Let $\Phi:\:G(\A)\rightarrow\C$ be an adelic function attached to $F$, where $G=\GSp(2n)$. If $F$ is an eigenform under all Hecke operators, then $\Phi$ generates an irreducible, automorphic representation of $G(\A)$.
\end{corollary}
{\bf Proof:} We apply Corollary \ref{secondglobaltheorem} with $G'=\SSp(2n,\R)^\pm$, the group of all $g\in\GSp(2n,\R)$ with multiplier $\pm1$. The identity component $G'_0=\SSp(2n,\R)$ has a maximal compact subgroup $K'_0$ isomorphic to $U(n)$. By our remarks above, condition i) of Corollary \ref{secondglobaltheorem} is satisfied. Considering weights, it is easy to see that conditions ii) and iii) in Corollary \ref{secondglobaltheorem} are satisfied. Condition iv) is equivalent to $F$ being holomorphic; see \cite{AS}. Condition v) is satisfied since $F$ is a Hecke eigenform.\qed

\begin{corollary}\label{secondglobaltheoremcor2}
 Let $F$ be a vector-valued holomorphic Siegel modular cusp form of degree $2$ and full level. Let $\Phi:\:G(\A)\rightarrow\C$ be an adelic function attached to $F$, where $G=\GSp(4)$. If $F$ is an eigenform under almost all Hecke operators, then $\Phi$ generates an irreducible, automorphic representation of $G(\A)$.
\end{corollary}
{\bf Proof:} By the argument in \cite{Ve}, $F$ is automatically an eigenform for \emph{all} Hecke operators. Hence the result follows from Corollary \ref{secondglobaltheoremcor1}.\qed

\end{document}